\documentclass[reqno]{amsart}

\usepackage{graphicx}
\usepackage{color}
\newtheorem{theorem}{Theorem}
 \newtheorem{lemma}{Lemma}
\newtheorem{corollary}{Corollary}

\UseRawInputEncoding
\usepackage[russian, english]{babel}
 \usepackage{setspace}

\begin{document}


\title[ 
Properties of  moments of density ]{
Properties of  moments of density for nonlocal mean field game equations
with a quadratic  cost 
}

\author{ Olga S. Rozanova*, Mikhail V. Inyakin}

\address{ Mathematics and Mechanics Department, Lomonosov Moscow State University, Leninskie Gory,
Moscow, 119991, Russian Federation, rozanova@mech.math.msu.su}



\begin{abstract}
We consider mean field game equations with an underlying  jump-diffusion process $X_t$ for the case of a quadratic cost function and show that
the expectation  and variance of $X_t$ obey  second-order ordinary differential equations with coefficients depending on the parameters of the cost function. Moreover, for the case of pure diffusion, the characteristic function and the fundamental solution of the equation for the probability density can be expressed in terms of the expectation ${\mathbb E}$ and the variance ${\mathbb V}$ of the process $X_t$, so that the moments of any order depend only on ${\mathbb E}$ and ${\mathbb V}$.
\end{abstract}

\maketitle



Mean field game theory (MFG) has been intensively developed in the last decade in terms of various mathematical disciplines such as partial differential equations, control theory, probability theory and numerical methods. It also has numerous practical applications in all areas where the behavior of a large number of agents is critical, e.g. in economics, mobile network design, logistics, etc. \cite{LL}, \cite{GLL}, \cite{Gomes}, \cite{Caines}.

In this paper we consider the forward-backward MFG  problem, which
analytically reduces to the solution of
the
initial-terminal problem for
coupled nonlocal Hamilton-Jacobi-Bellman (HJB) and Kolmogorov-Feller-Fokker-Planck (KFFP) equations for the density $m(t,x)$ and the value function $\Phi(t,x)$ \cite{LL}, \cite{Oksendal}, \cite{Oksendal1}:
\begin{eqnarray}
        & -\partial_t \Phi + \frac{1}{2} (|\nabla \Phi|)^2 -\frac{\delta^2}{2} \Delta \Phi - \lambda \left( \int\limits_{\mathbb R} \Phi(t,x+z) p(z)  dz - \Phi(t,x)  \right)   = g(t,x,m), \label{HJB}\\
        & \partial_t m - {\rm div} (m\nabla \Phi) - \frac{\delta^2}{2} \Delta m - \lambda \left( \int\limits_{\mathbb R} m(t,x-z) p(z)  dz - m(t,x)  \right)= 0, \label{KFP}\\
        & m(0,x) = m_0(x),\qquad \Phi(x,T) = K(x).\label{BC_KFP}
      \end{eqnarray}
      Here  $x\in {\mathbb R}^n$, $t\in [0, T]$, $0<T<\infty$. Let us stress that we do not discuss here the existence and smoothness of the solution and only assume that the initial and terminal conditions, as well as the right-hand side $g$ are compatible with the existence of moments of density of any order.

      In our consideration $m(t,x)\ge 0$ is the density of a stochastic process  $X$  with  dynamics given by
\begin{equation}\label{X}
dX_s = \alpha_s\, d s + \delta\, d W_s + \lambda d \Gamma_s, \quad X_0=x_0,
\end{equation}
 $x_0\in \mathbb R^n$  is a point in the space of states,
  $0\le s\le T$,
$W_s$ is a standard vectorial Brownian motion, $\Gamma_s$ is
 a purely discontinuous vector process with  jump size given by an independent random
variable,
 $\delta\ge 0$,  $\lambda\ge 0$
 are constants,
 $\int\limits_{{\mathbb R}^n}m(t,x)=1$.
 We denote $p(z)$  a probability density of jumps, $\int\limits_{{\mathbb R}^n}\,p(z) dz=1$, and assume
 $\mathcal M_1=\int\limits_{{\mathbb R}^n}\,z p(z) dz<\infty$, $\mathcal M_2=\int\limits_{{\mathbb R}^n}\,|z|^2 p(z) dz<\infty$.

Heuristically, the  problem is to minimize over all the progressively measurable admissible 
controls $\alpha_s\in L_{\mathcal F}^2(0,T;{\mathbb R}^n)$  the  cost function
$$ {\mathbb E}\Bigg[
\int_0^T \left( \frac{|\alpha_s|^2}{2} + g(s,X_s,m) \right)\,ds + K(X_T) \Bigg], $$
 $g:{\mathbb
R}_+\times {{\mathbb R}^n}\times  {\mathbb
R}_+\to \mathbb R$ and $K:{\mathbb R}^n\to
{\mathbb R}$  are prescribed continuous functions,
 the process  $X_s$ obeys (\ref{X}), the expectation is taken with respect to the filtration $\mathcal F$ generating by the jump-diffusion process.


 The main interest is the study of density evolution in the process of control.
      We are going to show that despite the impossibility of obtaining an explicit form of the density, for a particular case the behavior of its moments is relatively simple and often amenable to analytical study.

\bigskip

1. The main assumption that allows us to prove our main result is as follows:
\begin{equation}\label{g} g = a(t) |x|^2 + b(t)\cdot x +
c(t), 
\end{equation}
where $a,$ $ b=(b_1,\dots,b_n),$ $ c$  are smooth functions ($b$ is a vector) that can depend on the expectation ${\mathbb E}(X_s)$ and variance  ${\mathbb V}(X_s)$, and in this sense the control depends on the density,
and
\begin{equation}\label{K}
K(x)={ A}_T |x|^2+{ B}_T \cdot x +{ C}_T,
\end{equation}
with constant $A_T$, ${ B}_T$, $C_T$.

Similar assumptions about the cost function were made, for example, in \cite{Guean2009}, \cite{Bensoussan}, \cite{Yong}, \cite{Fatone} for the case of pure diffusion, where it was shown that the solution can be found by solving the matrix Riccati equation.
Our results are also based on the fact that the solution can be partially expressed in terms of the solution of the Riccati equation.

We are going to show that, with this special choice of control, the first two density moments can be found as solutions to second-order ODEs, and the characteristic function and fundamental solution are expressed in terms of them (and, perhaps, in terms of other parameters of the problem). Therefore, for any initial density, its full behavior can be restored.

\bigskip

2. First we find a solution to the HJB equation.
\begin{lemma}
 Assume that  $g$ has the form (\ref{g}). Then (\ref{HJB}) has a solution of the form
 \begin{equation}
\label{f2} \Phi = { A}(t)|x|^2 + { B}(t)\cdot x + { C}(t). \quad
\end{equation}
 The coefficients $A(t)$, $B(t)$, $C(t)$ are uniquely defined by the terminal conditions  (\ref{K}).
 \end{lemma}

\proof We substitute (\ref{f2}) to  (\ref{HJB}) and get a polynomial of the second order with respect to $x$. Then we combine the coefficients at $x^2$, $x$ and $1$ and obtain a system of ODEs for $A, B, C$, subject to the terminal conditions,
\begin{eqnarray}\label{AB}
        &&\dot{A} =2 A^2 - {a},  \quad
        \dot{B}-2 AB = 2 {\lambda} \mathcal M_1 \,
        A -{b},\\ && \dot{C} =-\frac12 B^2+\delta^2 A+\lambda (\mathcal M_2 A+\mathcal M_1 B), \nonumber \\
        && {A(T) = A_T}, \quad
        B(T) = B_T, \quad
        C(T) = C_T.\nonumber
         \end{eqnarray}
    The equation for $A$ splits from the rest of system and defines the entire dynamics. The equation for $B$ is linear and contain $A$ in coefficients.
$\Box$
\bigskip

It what follows we assume that $a(t)$ and $b(t)$ are such that
\begin{equation}\label{A_int}
e^{-2 \int\limits_0^T A(\tau) d\tau} <\infty, \qquad \int_{0}^{T}e^{-2\int\limits_\eta^T A(\tau) d\tau} B(\eta) \,d\eta\, <\infty.
\end{equation}


3. Now we find the {fundamental solution of the KFFP equation}. If $\Phi(t,x)$ is known, (\ref{KFP}) takes a more specific form,
\begin{eqnarray}\label{m1}
    \partial_t m - {( 2A \,x + B)} \cdot \nabla m  - {2A} m - \frac{\delta^2}{2} \Delta m - \lambda
    \left(   \int_{-\infty}^{\infty} m(t,x-z) p(z) dz - m  \right) = 0,
\end{eqnarray}
with initial condition
\begin{equation}\label{FS}
m_0(x) = \delta(x-y).
\end{equation}
We denote the solution to (\ref{m1}), (\ref{FS})  as $\mathcal G(t,x,y)$.

Applying the normalized inverse Fourier transform $x\to \omega$, we obtain the Cauchy problem for
the characteristic function
 $\psi=\psi(t,\omega,y)$:
\begin{eqnarray}
 \partial_t \hat{m} + 2{A(t)}\omega \partial_\omega \hat{m} +  \left(\frac{\delta^2}{2}  \omega^2 + i B(t)\omega  -   (\hat{p}(\omega)-1) \lambda \right) \hat{m} = 0,
   \quad \hat{m} (0,\omega) = e^{iy\cdot w},\label{meq}
\end{eqnarray}
where it is assumed that $ A $ and $ B $ are known from the previous step.
The solution to (\ref{meq}) can be explicitly found, namely,
\begin{eqnarray}\label{mhat}
 \psi (t,\omega,y)=
    {\exp}\left[-\int_{0}^{t}   \left(\frac{\delta^2}{2}
    {\mathcal R}^2 + i B(\eta)  {\mathcal R}  -   (\hat p({\mathcal R})-1)\lambda\right) d\eta+iy{\mathcal R}(t,0,\omega)\right],\\\nonumber
    \quad   {\mathcal R}= {\mathcal R}(t,\eta,\omega)= \omega e^{2 R(t,\eta)}, \quad R(t,\eta)=-\int\limits_\eta^t A(\tau) d\tau.
   \end{eqnarray}


\bigskip

4. We find the expectation and variance in terms of $A$ and $B$.

\begin{lemma}
Let  $X_t$  be a random process  with an initial probability density $m_0(x)$, $x_i^k \,m_0(x) \in L_1({\mathbb R)}$, $k=1,2$, $i=1,\dots,n$. The expectation and variance of $X_t$ can be found as
  \begin{equation}\label{T1E}
     \mathbb{E}(X_t)=\int_{0}^{t}e^{-2\int\limits_\eta^t A(\tau) d\tau} (B(\eta)+\lambda \mathcal M_1) \,d\eta + \mathbb{E}(X_0)\,e^{-2\int\limits_0^t A(\tau) d\tau},
\end{equation}
\begin{equation}\label{T1V}
    \mathbb{V}(X_t) =(\delta^2+\lambda \mathcal M_2)\,\int_{0}^{t}e^{-4\int\limits_\eta^t A(\tau) d\tau}  \,d\eta + \mathbb{V}(X_0)\,e^{- 4\int\limits_0^t A(\tau) d\tau},
\end{equation}
where $A(t), B(t)$ is a  solution to (\ref{AB}),   $t\in [0,T].$
\end{lemma}
  \proof
We denote as $m(t,x)$  the solution of \eqref{m1} subject to the initial data  $m_0(x)$. By means of standard computations we get from  (\ref{mhat})
 \begin{eqnarray*}\label{E}
 &&   \mathbb{E}(X_t) = -i \frac{\partial \hat{ m}(t,w)}{\partial \omega}\Big|_{w=0}  =
   - i \int\limits_{{\mathbb R}^n}\frac{\partial \psi (t,w,y)}{\partial \omega} \Big|_{w=0} m_0(y)\, dy=\\
 && e^{-2\int\limits_0^t A(\tau) d\tau}\,\int\limits_{{\mathbb R}^n} y m_0(y)\, dy+
   \int_{0}^{t}e^{2 R(t,\eta)} \left( B(\eta)+{\lambda}{\mathcal M}_1\right) \,d\eta,
\end{eqnarray*}
\begin{eqnarray*}\label{V}
 &&   \mathbb{V}(X_t) =  \mathbb{E}(X^2_t)- (\mathbb{E}(X_t))^2 =
    -\int\limits_{{\mathbb R}^n}\frac{\partial^2 \psi(t,w,y)}{\partial \omega^2} \Big|_{w=0} m_0(y)\, dy- (\mathbb{E}(X_t))^2=\\
 && e^{-4\int\limits_0^t A(\tau) d\tau}\, ( \int\limits_{{\mathbb R}^n} y^2 m_0(y)\, dy- (\mathbb{E}(X_t))^2)\,+
    (\delta^2+\lambda \mathcal M_2)\,\int_{0}^{t}e^{-4\int\limits_\eta^t A(\tau) d\tau}  \,d\eta.
    \end{eqnarray*}
Here we take into account that
$m(t,x)=\int\limits_{{\mathbb R}^n} {\mathcal G} (t,x,y)\, m_0(y) \, dy$ and
 $\int\limits_{{\mathbb R}^n} m_0(y)\, dy=1$. $\Box$

\bigskip

5. Now we express  $\psi$ in terms of $\mathbb E$ and $\mathbb V$.
\begin{theorem} Assume that condition \eqref{A_int} holds.
Then the characteristic function $\psi(t,w,y)$
 can be expressed in terms of $\mathbb E$, $\mathbb V$ and $A$ as follows:
 \begin{eqnarray*}\label{Ghat}
  &&\psi(t,w,y)=\\ 
  &&\exp\big[ -\frac12 |w|^2 ({\mathbb V} (X_t)-{\mathbb V} (X_0) {\mathbb K}^2(t))  + i w \cdot ( {\mathbb E} (X_t)-{\mathbb E} (X_0) {\mathbb K}(t)) + \lambda \mathcal Q (t, w) +i y\cdot w \,{\mathbb K}(t) \big]
  \end{eqnarray*}
  where $\mathbb K(t)= e^{-2 \int\limits_0^t A(\tau) d\tau}$,
   $Q (t, w)=\int\limits_0^t (\hat p({\mathcal R})-1- i  \mathcal M_1 \cdot {\mathcal R}) d\eta, \quad$
   ${\mathcal R}(t,\eta,\omega)= \omega e^{-2 \int\limits_\eta^t A(\tau) d\tau}$ ($\mathcal R$ and $ M_1$ are vectors).

   Moreover, $\mathbb K$ is a solution of the boundary problem
    \begin{eqnarray}\label{K}
  \ddot {\mathbb K} + 2 a \mathbb K=0, \quad \mathbb K(0)=1, \quad \mathbb K(T)=e^{-2 \int\limits_0^T A(\tau) d\tau}.
\end{eqnarray}

  For  $\lambda=0$ and $y=0$ the characteristic function depends on  $\mathbb E$, $\mathbb V$ only, therefore
 every moment of  $X_t$  can be expressed in terms of its expectation and variance.
  \end{theorem}

The {\it proof} follows from \eqref{mhat} and Lemma 2 directly, since
\begin{eqnarray*}\label{Ghat}
  {\psi}(t,w,y)= \exp\big[ -\frac12 |w|^2 \bar{\mathbb V} (X_t) + i w \cdot \bar{\mathbb E} (X_t) + \lambda \mathcal Q (t, w) + iy\cdot w \,e^{-2 \int\limits_0^t A(\tau) d\tau} \big],
  \end{eqnarray*}
  $\bar{\mathbb E}$ and $\bar{\mathbb V}$ are the expectation and variance of the process $X_t$, given as \eqref{X}, such that ${\mathbb E}(X_0)={\mathbb V}(X_0)=y=0$.
  Note that $A$ is defined by $a$ and the terminal condition, see \eqref{AB}.
  Property \eqref{K} follows from the first equation of \eqref{AB}.
$\Box$

\begin{corollary} In the case $\lambda=0$ the fundamental solution is
\begin{eqnarray*}\label{GG}
  {\mathcal G}(t,x,y)=\frac{1}{(\sqrt{2 \pi \bar{\mathbb V} (X_t)})^n}\,  \exp\big[ -\frac{|x-y \,e^{-2 \int\limits_0^t A(\tau) d\tau}-\bar{\mathbb E} (X_t)|^2}{2\bar {\mathbb V} (X_t)}\big].
  \end{eqnarray*}
\end{corollary}

\bigskip

6. Now we derive equations for the expectation and variance.
\begin{theorem} Assume that condition \eqref{A_int} holds.
 The expectation of  $\mathbb{E}(t)=\mathbb{E}(X_t)$ and the variance $\mathbb{V}(t)=\mathbb{V}(X_t)$ satisfy the following  boundary value problems:
   \begin{eqnarray}
\label{EE} &\mathbb{E}^{''}(t) + 2a(t)\mathbb{E}(t)=-b(t), \\\quad &\mathbb{E}(0)=\mathbb{E}(X_0), \quad
    \mathbb{E}(T) =\int_{0}^{T}e^{-2\int\limits_\eta^T A(\tau) d\tau} \left(B(\eta)+{\lambda}{\mathcal M}_1\right) \,d\eta + \mathbb{E}(X_0)\,e^{-2\int\limits_0^T A(\tau) d\tau},\nonumber
\end{eqnarray}
  \begin{eqnarray}
\label{Var} &
\mathbb{V}^{''}(t) + 4 a(t) \mathbb{V}(t)-\displaystyle\frac{(\mathbb{V}^{'}(t))^2 -K^2}{2 \mathbb{V}(t)}=0,
\quad K=\delta^2+\lambda {\mathcal M}_2,
 \\\quad &\mathbb{V}(0)=\mathbb{V}(X_0), \quad
    \mathbb{V}(T) =(\delta^2+\lambda \mathcal M_2)\,\int_{0}^{T}e^{-4\int\limits_\eta^T A(\tau) d\tau} \,d\eta + \mathbb{V}(X_0)\,e^{-4\int\limits_0^T A(\tau) d\tau},\quad t\in [0,T].\nonumber
\end{eqnarray}
\end{theorem}
\proof The statement follows from
\eqref{T1E} and \eqref{T1V} by direct computation. $\Box$

Note that  equation \eqref{EE} was obtained in \cite{RM}.

\bigskip

7. So far, we have not indicated the dependence of $a$ and $b$ on $\mathbb E$ and $\mathbb V$. We can use various hypotheses about this dependence and obtain system \eqref{EE}, \eqref{Var}, which is non-linear and, generally speaking, cannot be solved analytically. However, in several cases, the solution can be found explicitly. Corollary \ref{Cor2} \cite{RM}  gives one of these possibilities.

\begin{corollary}\label{Cor2}
If $a(t)=\rm const$ and $b(t)= b_0+b_1\mathbb E(t)+b_2\mathbb E'(t) $, $b_0, b_1, b_2=\rm const$,  then  (\ref{EE}) transforms into a linear second order ODE
  \begin{equation*}
\mathbb{E}^{''} + b_2\mathbb{E}'+(2a+b_1)\mathbb{E}=-b_0,
\end{equation*}
which can be explicitly solved.
\end{corollary}


   However, the simplest case of an explicit solution is  constants $a$ and $b$. Namely,

i) For $a>0$
\begin{eqnarray}
&&\mathbb{E}(t)=C_1 \sin \sqrt {2 t} +C_2 \cos \sqrt {2a} t-\frac{b}{2a},\label{aEp}\\
&&\mathbb{V}(t)=\pm \frac{1}{a}\,\sqrt{a \left((C_1^2+C_2^2) a+\frac{K^2}{8}\right)}\label{aVp}
+C_1 \sin(2\sqrt{2 a t})+ C_2 \cos(2\sqrt{2 a t}),
\end{eqnarray}
the sign is chosen for reasons of non-negativity $\mathbb{V}$.

In this case
\begin{eqnarray}
  A(t) &=& \sqrt{\frac{a}{2}} \,\tan \theta(t,T),\qquad \theta(t,T)=\arctan \sqrt{\frac{2}{a}} A_T+\sqrt{2a}(T-t),\nonumber\\
  B(t) &=& -\frac{b}{\sqrt{2a}}\,\tan \theta(t,T)+\frac{bA_T+aB_T}{\sqrt{a(a+2A_T^2)}\cos \theta(t,T)}.\nonumber
  \end{eqnarray}
Note that  the solution of \eqref{AB} may have singularities inside $(0, T)$, nevertheless condition \eqref{A_int} holds,
so the singularities of $A$ do not interfere with the calculation of the boundary conditions.

It is possible to show (see \cite{RN}), that
$\mathbb{E}(t)$ is differentiable in the points  $t$, where $A$ and $B$ fail to exist. In turn,
$\mathbb{V}(t)$ is continuous but not differentiable when $\mathbb{V}(t)=0$.

\bigskip

ii) For $a<0$
\begin{eqnarray}
&\mathbb{E}(t)=C_1 \sinh \sqrt {2 t} +C_2 \cosh \sqrt {2a} t-\frac{b}{2a},\label{aEm}\\
&\mathbb{V}(t)=
C_1+\frac{8 a C_1^2+K^2}{32 a C_2} \exp\left(2\sqrt{-2a}t\right)+ C_2\exp\left(-2\sqrt{-2a}t\right).\label{aVm}
\end{eqnarray}
Functions $A$ and $B$, necessary for boundary conditions, can be found from \eqref{AB}, the problem of singularities does not arise.

\bigskip

iii) For $a=0$ both $\mathbb{E}(t)$ and $\mathbb{V}(t)$ are polynomials of the second order with respect to $t$.

\bigskip

8.
Choosing $g$ in the form (\ref{g}) seems artificial. However, it appears naturally in a number of problems related to behavioral economics, e.g. \cite{TradeCrowding}, \cite{RN}, \cite{Trusov}. Note that in \cite{TradeCrowding} a second-order equation was obtained for the expectation under different assumptions. The MPG equations with nonlocal terms were considered in~\cite{Graber}.

It is well known that in the case of pure diffusion and a quadratic cost function, one can find an explicit solution for the density in the form of a Gaussian function \cite{Guean2009}.  In \cite{RN} an equation for the position of the maximum of density was obtained, which coincides for this case with the expectation.  For the case of a jump diffusion, these two characteristics of the density function do not coincide.

\bigskip

9. Since $\mathbb{E}(t)$ and $\mathbb{V}(t)$ are {observable} for real processes and can be estimated by statistical methods, we can propose a procedure for restoring the parameters of the cost functional under the assumption that the penalty term $g $ is quadratic. For example, assuming constant $a$ and $b$, one can use the formulas \eqref{aEp}, \eqref{aVp}
or \eqref{aEm}, \eqref{aVm} depending on the oscillatory or non-oscillatory character of $\mathbb{E}(t)$ and $\mathbb{V}(t)$ in time.
We can refer to \cite{Hommes}, \cite{Agliari} and their references to show that both types of behavior are possible.

\bigskip

10. Acknowledgments. Supported by the Moscow Center for
Fundamental and Applied Mathematics under the agreement
№075-15-2019-1621.

\end{document}